# Time-variant Nonlinear Participation Factors Considering Resonances in Power Systems

Tianwei Xia, *Student Member, IEEE*, Kai Sun, *Senior Member, IEEE*
Department of Electrical Engineering and Computer Science, University of Tennessee, Knoxville, TN, USA
txia4@vols.utk.edu, kaisun@utk.edu

*Abstract*— The participation factor (PF), as an important modal property for small-signal stability, evaluates the linkage between a state variable and a mode. Applying the normal form theory, a nonlinear PF can be defined to evaluate the participation of a state variable into modal dynamics following a large disturbance, that gives considerations to resonances and nonlinearities up to a desired order. However, existing nonlinear PFs are inconsistent with the conventional linear PF when nonlinear dynamics following a large disturbance attenuate and linear modal dynamics become dominating. This paper proposes a time-variant nonlinear PF by introducing a time decaying factor and the definition of a nonlinear mode. The new PFs consider modes of resonances and their values naturally transition to a linear PF when the system state becomes close to its equilibrium. The case study on a two-area four-generator system shows that the new PF can correctly rank generators by their participations in natural and resonance modes of nonlinear oscillation subject to a large disturbance.

## I. INTRODUCTION

Sustained power oscillations, especially those following a large disturbance, are threats to the stability and dynamic performance of a power system. In order to mitigate an oscillation to improve rotor angle stability of generators, control actions, e.g., damping control by power system stabilizers, can be taken at key generators that participate most in the oscillation mode. The participation factor (PF) is a useful index to evaluate the two-way linkage between each generator (or a state variable) and an oscillation mode (characterized by a pair of conjugate eigenvalues). It is defined as the product of corresponding elements on the mode shape and mode composition, and hence it considers both modal controllability and observability [1]. As a modal property, PFs are widely used in power system stability studies, such as PSS design [2], SVC placement [3].

A basic assumption with the PF is that the system model is linear, so it is usually used for linear system studies or small-signal stability analysis of a nonlinear system subject to a small disturbance. To extend the application of the conventional linear PF to nonlinear systems under a large disturbance, a nonlinear PF (NPF) is proposed in [3-7] based on the normal form theory, which transforms a nonlinear system to a linear system through a nonlinear coordinate transformation [3]. Paper [5] in 2005 summarizes the feasibility, contribution and application of the nonlinear PFs in detail.

In the past decade, PFs have been redefined or extended from different aspects, such as the extended PFs [8]-[10] and modal PFs [11]. However, many of these PFs are neither consistent nor unified in one framework. For instance, resonance modes are often considered by existing nonlinear PFs but are ignored by the conventional linear PF. Thus, when studying the process of damping an oscillation caused by a large disturbance, one often needs to switch from a nonlinear PF to the linear PF when judging the oscillation amplitude to be small enough. Thus, a discontinuity between nonlinear and linear PFs exists.

For a unified framework of PFs on oscillations following small and large disturbances, this paper proposes a new time-variant nonlinear PFs by introducing a time decaying factor and the definition of a nonlinear mode based on convolution of the spectrum and a normal distribution in the frequency domain. The new PF considers resonance modes, and naturally transitions to a linear PF when the system state becomes close to its equilibrium. Similar to the application based on nonlinear PF [6,7], this new time-variant PF as a kind of generalization of linear and nonlinear PFs can provide a reference for placement and configuration of PSS (Power System Stabilizers) for damping control under small and large disturbances.

The rest of the paper is organized as follows: Section II reviews the definitions of the linear PF, nonlinear PFs and other related PFs. Section III discusses the problems with existing nonlinear PFs in detail. Section IV introduces the proposed time-variant nonlinear PF. Besides, the convolution is used to include influences from resonance modes. The case study on Kundur's two-area system is presented in section V to show that new PFs can correctly rank generators by their participations in natural and resonance modes of nonlinear oscillation subject to a large disturbance. At last, the paper concludes in section VI.

## II. BACKGROUND OF DIFFERENT PARTICIPATION FACTORS

### A. *The Conventional Linear PFs* [1]

A power system can be modeled by a set of nonlinear differential equations:

$$\dot{\mathbf{x}} = f(\mathbf{x}, \mathbf{u}) \quad (1)$$

Linearize the model at its equilibrium point to generate:

$$\dot{\mathbf{x}} = \mathbf{A}\mathbf{x} \quad (2)$$

where $\mathbf{A}$ is the state matrix with the $i$-th eigenvalue or mode of

This work was supported in part by the ERC Program of the NSF and U.S. DOE under grant EEC-1041877 and in part by the NSF grant ECCS-1553863.

state denoted by $\lambda_i$. Its left eigenvector $\boldsymbol{\psi}_i=[\psi_{i1}, \ldots \psi_{in}]$ and right eigenvector are denoted by $\boldsymbol{\phi}_i=[\phi_{1i}, \ldots \phi_{ni}]^T$, respectively. They are actually the mode composition and mode shape of the $i$-th mode of $\lambda_i$. The conventional linear PF of the $k$-th state in the $i$-th mode is dimensionless as defined by

$$p_{ki} \triangleq \phi_{ki}\psi_{ik} \qquad (3)$$

which considers both the mode shape and mode composition, or in other words, both the activity of the state variable in mode and its contribution to the mode.

From initial state $\mathbf{x}_0$, responses of each state variable are

$$x_k(t) = \sum_{i=1}^{n}(\boldsymbol{\psi}_i\mathbf{x}_0)\phi_{ki}e^{\lambda_i t} = \sum_{i=1}^{n}B_{ki}e^{\lambda_i t} \quad (k=1,\ldots,n) \qquad (4)$$

where $B_{ki}=\boldsymbol{\psi}_i\mathbf{x}_0\phi_{ki}$ is named the contribution factor, which evaluates how much mode $i$ is excited in the $k$-th state variable. Based on (4), here is an interpretation of PFs defined in (3):

*The participation factor represents the size of the modal dynamics with a state variable when only that state variable is perturbed [4].*

The PF can be viewed as the contribution factor $B_{ki}$ when

$$\mathbf{x}_0 = \mathbf{e}_k = \begin{bmatrix} 0 & \ldots & 1|_k & \ldots & 0 \end{bmatrix}. \qquad (5)$$

The response of the system for the $k$-th state variable will be

$$x_k(t) = \sum_{i=1}^{n}(\boldsymbol{\psi}_i\mathbf{e}_k)\phi_{ki}e^{\lambda_i t} = \sum_{i=1}^{n}p_{ki}e^{\lambda_i t} \qquad (6)$$

This interpretation of the PF is extended in [3-7] as well as in the rest of this paper.

*B. Normal Form Theory and Nonlinear PFs*

The normal form method is a powerful tool to analyze nonlinearities of desired orders, and its idea is to transform a nonlinear system to a linear system through a nonlinear coordinate transformation. Theoretically speaking, the normal form method can be conducted for any order, but the nonlinearity of the 2nd order is widely considered [5].

Apply Tayler expansion to (1) at its stable equilibrium and keep the 1st order and 2nd order terms as shown in

$$\dot{\mathbf{x}} = A\mathbf{x} + \frac{1}{2}\begin{bmatrix} \mathbf{x}^T\mathbf{H}^1\mathbf{x} \\ \vdots \\ \mathbf{x}^T\mathbf{H}^n\mathbf{x} \end{bmatrix} \qquad (7)$$

where $\mathbf{H}$ is the Hessian matrix. Apply transformation $\mathbf{x}=\Phi\mathbf{y}$. The system in $\mathbf{y}$-space becomes

$$\dot{y}_i = \lambda_i y_i + \sum_{p=1}^{n}\sum_{q=1}^{n}C_{pq}^i y_p y_q \qquad (8)$$

where $C_{pq}^i$ are the coefficients of 2nd order terms after the transformation. In order to obtain a linear system, a nonlinear coordinate transformation from $\mathbf{y}=\mathbf{h}(\mathbf{z})$ is introduced [2]:

$$y_i = z_i + \sum_{p=1}^{n}\sum_{q=1}^{n}h2_{pq}^i z_p z_q \qquad (9)$$

$$h2_{pq}^i = \frac{C_{pq}^i}{\lambda_p + \lambda_q - \lambda_i}$$

Plug (9) into (8) to eliminate all 2nd order terms, and then the system in $z$-space becomes a formally more linear system, whose nonlinearities only appear on terms of the 3rd order or higher. The expressions in $x$, $y$ and $z$ space are [3]

$$z_i(t) = z_{i0}e^{\lambda_i t}$$

$$y_i(t) = z_{i0}e^{\lambda_i t} + \sum_{p=1}^{n}\sum_{q=1}^{n}h2_{pq}^i z_{p0}z_{q0}e^{(\lambda_p+\lambda_q)t} \qquad (10)$$

$$x_k(t) = \sum_{i=1}^{n}\phi_{ki}z_{i0}e^{\lambda_i t} + \sum_{i=1}^{n}\phi_{ki}[\sum_{p=1}^{n}\sum_{q=1}^{n}h2_{pq}^i z_{p0}z_{q0}e^{(\lambda_p+\lambda_q)t}]$$

By the interpretation in (6), let $\mathbf{x}_0 = \mathbf{e}_k$, and the initial state $z_{i0}$ can be approximated by

$$z_{i0} = \psi_{ik} - \sum_{p=1}^{n}\sum_{q=p}^{n}h2_{pq}^i \psi_{pk}\psi_{qk} \qquad (11)$$

Substitute (11) into (10). The system response becomes [4]:

$$x_k(t) = \sum_{i=1}^{n}p_{2ki}e^{\lambda_i t} + \sum_{p=1}^{n}\sum_{q=p}^{n}p_{2kpq}e^{(\lambda_p+\lambda_q)t}$$

$$p_{2ki} = \phi_{ki}(\psi_{ik} + \psi_{2ikk}) = p_{ki} + p_{2kiNL} \qquad (12)$$

$$p_{2kpq} = \phi_{2kpq}(\psi_{pk} + \psi_{2pkk})(\psi_{qk} + \psi_{2qkk})$$

where

$$\psi_{2mkk} = -\sum_{p=1}^{n}\sum_{q=p}^{n}h2_{pq}^m \psi_{pk}\psi_{qk} \qquad \phi_{2kpq} = \sum_{i=1}^{n}h2_{pq}^i \phi_{ki} \qquad (13)$$

In (12), the $p_{2ij}$ is the nonlinear PF of state $k$ in mode $i$. Notice that $p_{ki}$ is the linear PF, so the nonlinear PF $p_{2ki}$ can be viewed as the sum of the linear PF and an extra term $p_{2ijNL}$. As for the $p_{2kpq}$, it is related to the participation of the state in the resonance mode $\lambda_p + \lambda_q$, and is usually ignored [5]. This paper will discuss this ignored resonance mode in detail in section III.

*C. Other Participation Factors*

In [8] and [9], an extended PF considering a set of initial states is proposed as

$$\hat{p}_{ki} = \underset{x_{k0}^{(l)} \in S_k}{avg} \frac{B_{ki}}{x_{k0}^{(l)}} \qquad (14)$$

where $S_k$ is the set of the initial states that the $k$-th state variable can take in those initial states, and "$(l)$" indicates the $l$-th value of the set $S_k$. It evaluates the average linear contribution over a set $S_k$ of the initial states.

This definition is extended to the nonlinear PF in [9] to be related to initial states and the excitation energy when $t = 0$.

III. CONNECTION BETWEEN LINEAR AND NONLINEAR PFS

The inconsistency between linear and nonlinear PFs is caused by three factors, which will be discussed as follows.

*A) On Scaling on the Initial State*

The linear PF defined in (3) based on model (2) is dimensionless and does not depend on the size of the initial state. Although the interpretation (6) involves the initial state, its effect can be eliminated after normalization. For instance, in (5), if the amplitude of the initial state is scaled by $\alpha_k$, i.e. (15), the response will also be scaled by $\alpha_k$ in (16).

$$\mathbf{x}_0 = \alpha_k\mathbf{e}_k = \begin{bmatrix} 0 & \ldots & \alpha_k|_k & \ldots & 0 \end{bmatrix} \qquad (15)$$

$$B_{ki} = p_{ki} = \alpha_k\phi_{ki}\psi_{ik} \qquad (16)$$

Ref [7] assumes the scaling factor $\alpha_k = \alpha$ for any $k$. This means that state variables are assumed to have the same size of

excitation. There is

$$p_{ki} = \alpha_k \phi_{ki} \psi_{ik} = \alpha \phi_{ki} \psi_{ik} \quad (17)$$

The scaling factor can be canceled after the normalization of all PFs for any $\alpha$. Thus, the influence from the initial state is eliminated.

However, for a nonlinear PF defined by, e.g., (18), such a scaling factor cannot be eliminated by normalization.

$$p_{2ki} = \phi_{ki}(\alpha \psi_{ik} + \psi_{2ikk}) = \alpha \phi_{ki} \psi_{ik} - \alpha^2 \sum_{p=1}^{n} \sum_{q=p}^{n} h2_{pq}^m \psi_{pk} \psi_{qk} \quad (18)$$

Notice that for generator $k$, a nonlinear PF changes with $\alpha$ after the normalization. When $\alpha \to 0$, $p_{2ki} \to p_{ki}$. Papers [3-7] simply assume $\alpha = 1$. This paper will address the influence from the scaling on the initial state.

*B) On Discontinuity Between Linear and Nonlinear PFs*

Both the linear PF and existing nonlinear PFs are computed based on the system model and are invariant in time. Thus, when studying the dynamic process of damping an oscillation caused by a large disturbance, we have to switch from a nonlinear PF to the linear PF when the oscillation amplitude becomes small enough. This results in a discontinuity between two PFs and a practical problem: Which PF to trust? An engineering solution could be to calculate both and choose the nonlinear or linear PF by whether the oscillation exhibits a strong nonlinearity or resonance.

The introduction of the scaling factor $\alpha$ seems to be able to address this problem by introducing adjusting weights depending on the value of $\alpha$ for linear and nonlinear PFs from the comparison of (17) and (18). A nonlinear PF is more emphasized than a linear PF with a larger $\alpha$ when, for example, the system state is far from the stable equilibrium. However, consider a scenario that the system state oscillates starting from the edge of the domain of attraction, and is gradually damped to reduce its amplitude, and be closer and closer to the stable equilibrium. Consistency between nonlinear and liner PFs will be important for designing a unified control strategy based on PFs. Thus, it is desired to satisfy (19) to make a nonlinear PF smooth transition to the linear PF. Thus, the time of the transitioning will have to be concerned.

$$\lim_{t \to \infty} p_{2ki} = p_{ki} \quad (19)$$

*C) On Resonance Modes*

In (12), $p_{2kpq}$ is considered the participation of state $x_k$ in a mode of resonance between $\lambda_p$ and $\lambda_q$, which is usually ignored in small-signal stability analysis and control even if it is sometimes evaluated as an additional index to explain resonance phenomena [5, 6]. Existing works focus on the linear or nonlinear PF for each linear mode $\lambda_k$. However, the study in [11] shows that resonance modes can be easily observed from measurements, especially under large disturbances and may impact the system stability.

### IV. PROPOSED TIME VARIANT NONLINEAR PFS ADDRESSING RESONANCES

To address the three factors in section III, two steps are taken. First, introduce a time decaying factor to address factors 1 and 2. Second, use convolution to define a nonlinear mode considering the influence from resonances. Finally, time-variant nonlinear PF (TNPF) is obtained.

*A. Time Decaying Factor*

Based on the interpretation (6) and discussion in III-A, the participation factor can be viewed as the excitation energy of the selected mode $i$ when only a certified generator $k$ is perturbed with a time decaying factor $\alpha_k e^{\lambda_i t}$ rather than a constant. Its value is calculated by

$$p_{2ki}(t) = (\alpha_k e^{\lambda_i t})\phi_{ki}\psi_{ik} - \phi_{ki}\sum_{p=1}^{n}\sum_{q=p}^{n}(\alpha_p \alpha_q e^{(\lambda_p+\lambda_q)t})h2_{pq}^m \psi_{pk}\psi_{qk} \quad (20)$$

Equivalently, its value can also be estimated from a set of initial states by using the method of the extended PF

$$\tilde{p}_{2ki}(t) = \underset{\alpha_k^{(l)} \in S_k}{avg} \frac{p_{2ki}(t)}{\alpha_k^{(l)}} \quad (21)$$

The set $S_k$ can be obtained from the history measurement data. From the discussion in III-A, the value of the new defined PF depends on the value of $\alpha_k$. Obviously, when $t = 0$ and $\alpha_k = 1$, the new defined PF equals the nonlinear PF in (12). Without better knowledge, $\alpha_k=1$ is advised and will be used in the case study of this paper.

*B. Resonance Modes and Nonlinear Modes*

To include the 2$^{nd}$ order resonance terms in (12), the concept of mode needs to be extended. In the conventional definition, eigenvalue $\lambda_k$ (together with its conjugate) corresponds to one linear mode; the resonance mode corresponds to $\lambda_p + \lambda_q$, resulting from the 2$^{nd}$ nonlinearity. For 2$^{nd}$ order resonance terms, it would be valuable to define a nonlinear mode to reflect the contribution of 2$^{nd}$ order resonances.

**Definition 1** (*Nonlinear Mode*) A nonlinear mode $\mu_j$ is the weighted sum of all linear modes and resonance modes whose frequencies belong to a given range of $(f_{j,lower}, f_{j,upper}]$.

$$\mu_j = \sum_i w_i p_{ki} e^{\lambda_i t} \quad (22)$$

for all $i$'s satisfying $\text{Im}(\lambda_i) \in (f_{j,lower}, f_{j,upper}]$

where $w_i$ is the weight, $f_{j,lower}$ and $f_{j,upper}$ are the lower and upper frequency limit, respectively. For a certain state variable, it is easy to plot the spectrum of nonlinear PF for linear modes and the PF for resonance modes. Then the spectrum can be divided into several frequency intervals. The sum of linear modes and resonance modes in each interval is viewed as one nonlinear mode $\mu_j$.

For the weights $w_i$, a relatively reasonable way is to use the normal distribution. It is known that a weighted moving average is a convolution. In other words, for a certain frequency $f_j = (f_{j,lower} + f_{j,upper})/2$, the nonlinear mode will be

$$\mu_j = \int_0^\infty N(f_j, \sigma^2) g(f) df$$

$$g(f) = \begin{cases} p_{2ki} e^{\lambda_i t} & f \in \{\text{Im}(\lambda_i)\} \\ p_{2kpq} e^{(\lambda_p + \lambda_q)t} & f \in \{\text{Im}(\lambda_p + \lambda_q)\} \\ 0 & \text{others} \end{cases} \quad (23)$$

$g(f)$ is the spectrum of linear modes and resonance modes for a certain state variable.

### C. Time-variant Nonlinear PF (TNPF)

Part *A)* introduces the time decaying factor to connect linear and nonlinear PFs. In part *B)*, the convolution operation is used to include the influence of resonance modes. To take both advantages of the two aspects, the time-variant PFs are defined as

$$p_2(t, f_{taregt}) = \int_0^\infty N(f_{taregt}, \sigma^2) p_g(f) df \quad (24)$$

$$p_g(f) = \begin{cases} p_{2ki}(t) & f \in \{\text{Im}(\lambda_i)\} \\ p_{2kpq}(t) & f \in \{\text{Im}(\lambda_p + \lambda_q)\} \\ 0 & \text{others} \end{cases}$$

where

$$p_{2ki}(t) = (\alpha_k e^{\lambda_i t})\phi_{ki}\psi_{ik} - \phi_{ki}\sum_{p=1}^n\sum_{q=p}^n (\alpha_p\alpha_q e^{(\lambda_p+\lambda_q)t}) h2_{pq}^m \psi_{pk}\psi_{qk} \quad (25)$$

$$p_{2kpq}(t) = (\alpha_p\alpha_q e^{\lambda_p t + \lambda_q t})\phi_{2kpq}(\psi_{pk} + \psi_{2pkk})(\psi_{qk} + \psi_{2qkk})$$

$P_g(f)$ presents a participation-factor spectrum for linear modes and resonance modes. $P_{2ki}(t)$ is the time-variant PF for the linear mode, and $P_{2kpq}(t)$ is the time-variant PF for the resonance mode.

Fig. 2 illustrates the spectrum of nonlinear PFs and TNPFs. For convenience, $t = 0$ here. Since the TNPFs can be viewed as the convolution of nonlinear PFs and the normal distribution, the standard deviation σ will be critical to the performance. It is related to the resolution of the measurement and control side. The detailed study of this is out of scope in this paper, and the case study shows that σ = 0.1 is a good selection.

In (25), $P_{2kpq}(t)$ contains 3rd order and 4th order decaying terms because of $\psi_{pkk}$ and $\psi_{qkk}$. For simplification, those terms are all approximated by the 2nd order decaying terms, since $\psi_{pkk}$ and $\psi_{pkk}$ are usually much smaller than $\psi_{pk}$ and $\psi_{qk}$. All PFs are summarized in TABLE I.

TABLE I. COMPARISON OF DIFFERENT PFs

| Types | Time Performance | System model | Mode | Formula |
|---|---|---|---|---|
| Linear PF | Constant | Linear | Linear | (3) |
| Nonlinear PF | Constant | Nonlinear | Linear | (12) |
| TNPF | Time-variant | Nonlinear | Nonlinear | (24) |

## V. CASE STUDY

Kundur's two-area system is used for the test, whose topology is shown as

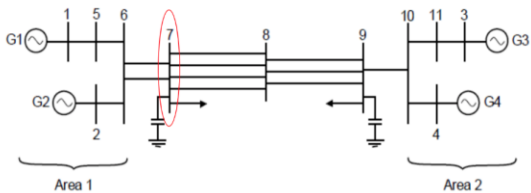

Fig. 1. Kundur's two-area system [1]

The linear and nonlinear PFs are shown in Table II. 1.11 Hz and 1.61 Hz modes are local modes, where the linear and nonlinear PFs are almost the same. 0.59 Hz mode is the inter-area mode, for which the ranking of 4 generators is different. Generator 3 is the second most important generator by linear PFs, while generator 4 is the second by nonlinear PFs. Since they are all time-invariant constants, without the TNPF, the system operator will have to choose between generators 3 and 4 for a control action.

TABLE II. LINEAR AND NONLINEAR PFs FOR TWO-AREA SYSTEM

| Gen. | 0.59 Hz | | 1.11 Hz | | 1.61 Hz | |
|---|---|---|---|---|---|---|
| | Linear PF | Nonlinear PF | Linear PF | Nonlinear PF | Linear PF | Nonlinear PF |
| 1 | 1.00 | 1.00 | $2.05\times10^{-3}$ | $2.06\times10^{-3}$ | 0.13 | 0.11 |
| 2 | 0.05 | 0.02 | $6.20\times10^{-3}$ | $5.97\times10^{-3}$ | 1.00 | 1.00 |
| 3 | 0.86 | 0.37 | 0.80 | 0.81 | $2.86\times10^{-3}$ | $2.76\times10^{-3}$ |
| 4 | 0.62 | 0.40 | 1.00 | 1.00 | $8.31\times10^{-3}$ | $7.69\times10^{-3}$ |

To show the performance of the TNPF, first select generator 1. Fig. 2 shows the spectrum about its PFs, nonlinear PFs and TNPF ($t = 0$). The blue triangles and black squares are the linear and nonlinear PFs for the linear modes. Notice that they share the same frequency. The red circles show the nonlinear PFs for the resonance mode. This system is carefully designed to make the 1.11 Hz mode resonant with the harmonic mode of 0.59 Hz mode. The red curve shows the convolution of the normal distribution.

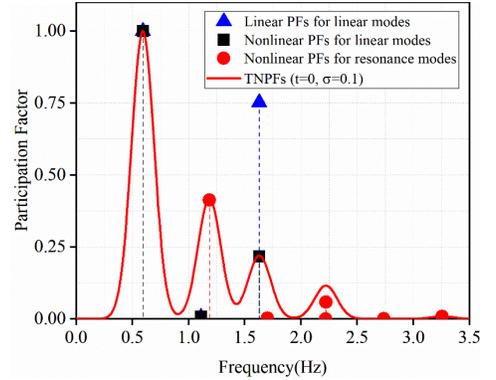

Fig. 2. Spectrum of generator 1

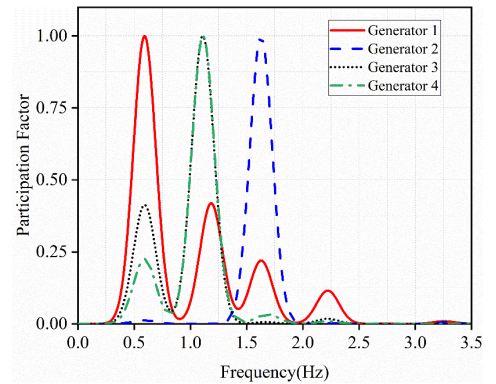

Fig. 3. Spectrum on all four generators

Next, the TNPFs (t = 0) for all generators are shown in Fig.3. It is obvious that four nonlinear modes are obtained through convolution. The nonlinear mode of 2.2 Hz is the resonance mode, and it comes from the harmonic terms of 1.11 Hz mode and the interaction between 1.61 Hz mode and 0.59 Hz mode.

On the 0.59 Hz inter-area mode, the difference between linear and nonlinear PFs in values results in a discontinuity. The proposed TNPF can address this problem. The result of TNPFs is shown in Fig. 4. The TNPF is the same as the nonlinear PF when $t = 0$. With the increasing of $t$, due to the decaying terms, the TNPF becomes closer and closer to the linear PF, which coincides with the fact that for a stable system, the system will become more and more linear after the fault. There is no resonance mode near 0.59 Hz, so the $P_{2kpd}(t)$ for 0.59 Hz is almost zero.

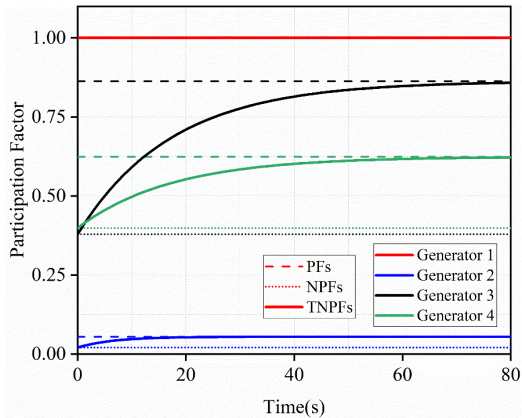

Fig. 4. The trajectories of PF, nonlinear PF and TNPF for 0.59 Hz mode

The 1.11 Hz mode has a different situation. For the 1.11 Hz linear mode, the linear and nonlinear PFs will share almost the same result (in TABLE II). However, for the 1.11 Hz nonlinear mode, the resonance mode will significantly influence the final result, which can be viewed from the resonance mode (1.1 Hz red circle) in Fig.2. The TNPFs of the 1.11 Hz nonlinear mode in the time domain are shown in Fig.5. The TNPF for generator 1 is much larger than the linear or nonlinear PF when $t = 0$ due to the resonance mode. With the time increasing, the resonance mode has a large damping term, so the value of TNPF decreases fast. The stable result will be the same as the linear PF.

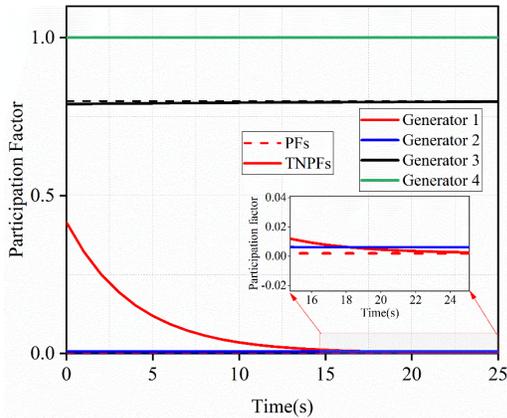

Fig. 5. The trajectories of PF and TNPF for 1.11 Hz mode

## VI. CONCLUSION

This paper has discussed the inconsistency between linear PFs and nonlinear PFs, which are caused by three factors. To address those factors, the TNPF is proposed. The time decaying factor in TNPF is used to allow a smooth transition from the NPF to the linear PF for an oscillating power system subject to a large disturbance and to address resonances. By introducing the convolution operation, nonlinear modes can be defined to calculate the TNPF. The case study on a two-area system has demonstrated the performances of the new PFs in ranking generators by their participations in power system oscillations energized by a large disturbance. Our future work will apply the proposed TNPF to power system stability monitoring and control under large disturbances when nonlinearities in power system oscillations cannot be ignored [12]-[15].

REFERENCES

[1] P. Kundur, *Power System Stability and Control*. Palo Alto, CA: McGraw-Hill, 1993, pp. 715–723

[2] L. Hassan, M. Moghavvemi, et al., "Optimization of power system stabilizers using participation factor and genetic algorithm," *International Journal of Electrical Power & Energy Systems*, vol. 55, no. 4, pp. 668-679,

[3] J. Zhang, J. Y. Wen, S. J. Cheng et al., "A Novel SVC Allocation Method for Power System Voltage Stability Enhancement by Normal Forms of Diffeomorphism," *IEEE Trans. Power Syst.*, vol. 22, no. 4, pp. 1819-1825, 2007.

[4] S. K. Starrett and A. A. Fouad, "Nonlinear measures of mode-machine participation," *IEEE Trans. Power Syst.*, vol. 13, no. 2, pp. 389–394, May 1998.

[5] J. J. Sanchez-Gasca, V. Vittal, et al., "Inclusion of higher order terms for small-signal (modal) analysis: committee report-task force on assessing the need to include higher order terms for small-signal (modal) analysis," *IEEE Trans. Power Syst.*, vol. 20, no. 4, pp. 1886-1904, Oct. 2005

[6] S. Liu, A. R. Messina, et al., "Assessing placement of controllers and nonlinear behavior using normal form analysis," *IEEE Trans. Power Syst.*, vol. 20, no. 3, pp. 1486-1495, Aug. 2005

[7] S. Liu, A. R. Messina, et al., "A Normal Form Analysis Approach to Siting Power System Stabilizers (PSSs) and Assessing Power System Nonlinear Behavior," *IEEE Trans. Power Syst.*, vol. 21, no. 4, pp. 1755-1762, Oct. 2006

[8] E. H. Abed, D. Lindsay, et al, "On participation factors for linear systems," *Automatica*, vol. 36, no. 10, pp. 1489-1496, Oct. 2000.

[9] B. Hamzi and E. Abed, "Local modal participation analysis of nonlinear systems using Poincaré linearization," *Nonlinear Dynamics*, vol. 99, no. 1, pp. 803-811, Jan. 2020.

[10] G. Tzounas, I. Dassios, F. Milano, "Modal Participation Factors of Algebraic Variables," *IEEE Trans. Power Syst.*, vol. 35, no. 1, pp. 742-750, Jan. 2020

[11] T. Xia, Z. Yu, et al., "Extended Prony Analysis on Power System Oscillation Under a Near-Resonance Condition," 2020 *IEEE Power & Energy Society General Meeting*, 2020, pp. 1-5.

[12] B. Wang, K. Sun, "Formulation and Characterization of Power System Electromechanical Oscillations," *IEEE Trans. Power Syst.*, vol. 61, no. 6, pp. 5082-5093, Nov. 2016

[13] X. Xu, W. Ju, B. Wang, K. Sun, "Real-time Damping Estimation on Nonlinear Electromechanical Oscillation," *IEEE Trans. Power Syst.*, vol. 36, No. 4, pp. 3142-3152, July 2021.

[14] X. Xu, K. Sun, "Direct Damping Feedback Control Using Power Electronics-Interfaced Resources," *IEEE Trans. Power Syst.*, In Press.

[15] Y. Zhu, C. Liu, K. Sun, et al, "Optimization of Battery Energy Storage to Improve Power System Oscillation Damping," *IEEE Trans. Sustainable Energy*, vol. 10, No. 3, pp. 1015-1024, July 2019.